\documentclass[11pt]{amsart}
\usepackage{amssymb}

\title{Questions on self maps of algebraic varieties}

\author{Najmuddin Fakhruddin}
 
\address{School of Mathematics, Tata Institue of Fundamental Research, Homi Bhabha Road,
Mumbai 400005, India}

\email{naf@math.tifr.res.in}

\newcommand{\te}{\otimes}
\newcommand{\bP}{\mathbb{P}}
\newcommand{\Q}{\mathbb{Q}}
\newcommand{\Z}{\mathbb{Z}}

\newcommand{\R}{\mathbb{R}}
\newcommand{\mc}{\mathcal}
\newcommand{\mb}{\mathbb}
\newcommand{\mf}{\mathbf}

\newcommand{\omx}{\omega_X}

\newtheorem{thm}{Theorem}[section]
\newtheorem{prop}[thm]{Proposition}
\newtheorem{conj}[thm]{Conjecture}
\newtheorem{cor}[thm]{Corollary}
\newtheorem{lem}[thm]{Lemma}

\theoremstyle{definition}
\newtheorem{ques}[thm]{Question}
\newtheorem{rem}[thm]{Remark}

\input{xypic}
\begin{document}

\maketitle 

\section{Introduction} \label{sec:intro}
In this note we shall consider some arithmetic and geometric questions
concerning projective varieties $X$ with a self-map $\phi$
and an ample  line bundle ${{L}}$ such that $\phi^*({{L}}) \cong {{L}}^{\te d}$
for some $d>1$, all defined over some field $k$. 
Such a situation is interesting arithmetically,
if $k$ is a number field, because
using a method of Tate (see for example \cite[Chapter 3]{serre-mw})
one can define a canonical height function $h_\phi:X(\bar{k}) \to \R$
which can be used to study the dynamics of $\phi$ on the set of rational points
of $X$. One is then naturally led to ask which varieties admit such self maps
and line bundles.

Basic examples of such varieties are projective spaces and
abelian varieties. The results of Section \ref{sec:1} show that projective spaces
are in a sense the ``universal'' such varieties, so several questions like the 
uniform boundedness of the torsion of abelian varieties over a number field
or the generalisation of Bogomolov's conjecture 
due to  Zhang \cite[Conjecture 2.5]{zhang-adelic}, reduce to questions
about self maps of projective spaces. (Unfortunately, the questions do not seem to
become any easier!)

In the next section we raise some questions on the arithmetic of self maps of 
projective spaces, among them a generalisation of the Morton-Silverman
Uniform Boundedness conjecture \cite{morton-silverman}. Most of these
seem beyond reach at the moment, however, we do show 
(Proposition \ref{prop:lang}) that
the weakest of these questions is implied by a conjecture of Lang, via the
work of Caporaso, Harris and Mazur \cite{chm}.

The remaining two sections are devoted to purely geometric questions.
In Section \ref{sec:3} we obtain a partial classification
(Theorem \ref{thm:kodzero}) of varieties
admitting self maps and line bundles as above. This suggests 
that all \'etale self maps ``come'' from maps of families of abelian
varieties. We also ask if
projective spaces can be characterised in terms of self maps.
Section \ref{sec:4} is devoted to showing that the self maps 
considered above have many periodic points; in fact the periodic points
always form a Zariski dense set.

\section{} \label{sec:1}

\begin{prop}
Let $X$ be a projective variety over an infinite field $k$, $\phi:X \to X$
a morphism
 and ${{L}}$  a
very ample line bundle on $X$ such that 
$\phi^*({{{L}}}) \cong {{{L}}}^{\te d}$ for some $d > 1$. 
Let $V = H^0(X, {{L}})$ and 
suppose the following conditions are satisfied:
\begin{enumerate}
\item The embedding $\iota:X \to \bP(V)$ induced by ${{L}}$
has the property that the maps $H^0(\bP(V),\mc{O}(n)) \stackrel{\iota^*}{\longrightarrow}
H^0(X, {{L}}^{\otimes n})$ are surjective for all $n \geq 0$. 
\item $\iota(X)$ is cut out set-theoretically by
homogenous forms of degree  $ \leq d$.
\end{enumerate}
Then there exists a morphism $\psi:\bP(V) \to \bP(V)$ such that $ \psi \circ \iota  = \iota 
\circ \phi$. \label{prop:embed}.
\end{prop}

\begin{proof}
Let $g =\dim(X)$, let $s_0,s_1,\ldots,s_n$ be a basis of $V$ and assume that 
$s_0,s_1,\ldots,s_g$ have no
common zeros on $X$. By (1), it follows that after choosing an isomorphism 
$\phi^*({{{L}}}) \cong {{{L}}}^{\te d}$, each $\phi^*(s_j)$, $0 \leq j \leq n$,
 can be written as a homogenous polynomial of degree $d$ (not necessarily uniquely) in the 
$s_i$'s. Any choice of such polynomials $f_j$, $0 \leq j \leq n$, gives rise to a rational
 map  $\psi :  \bP(V) \to \bP(V)$ defined by $s_j \mapsto f_j$, such that  
$ \psi \circ \iota  = \iota \circ \phi$.  We shall show that if the $f_j$'s are chosen generally
then $\psi$ is actually a morphism.

Let $ D(f_j)$ denote the hypersurface in $\bP(V)$ which is the zero set of $f_j$.
Using induction we will prove that 
 the intersection of all the $ D(f_j)$'s is empty.
Firstly, each component of the intersection of the first
$g+1$ $ D(f_i)$'s is of codimension $ g+1$, since otherwise it would intersect $\iota(X)$,
and by choice of the $s_j$'s this cannot happen. 
Now we prove that we can 
choose the remaining $f_i$'s in such a way so that the codimension of the
intersection goes down by $1$ at each step, this clearly suffices to complete the
proof.

Suppose $g \leq k < n$, we have chosen $f_i$ for 
for all $ i \leq k $ and we want to choose $f_{k+1}$. Let $h_1$ be any polynomial
of degree $d$ which restricts to  $\phi^*(s_{k+1})$ on $A$. If $D(h_1)$ does not 
contain any of the components of the intersection of the $D(f_j)$'s, $j \leq k$,
then let  $f_{k+1} = h_1$. Otherwise, by (2) we can choose a polynomial $h_2$ of degree $d$
in the ideal of $X$ such that $D(h_2)$ does not contain any of the above components. Consider 
the family of hypersurfaces $D(ah_1 + bh_2)$, $[a,b] \in
{\bP}^1$. When $a=0$, by construction the corresponding hypersurface does not contain any of 
the components. Since  this is an open condition and $k$ is infinite,  
there exists $b \in k$ such that the 
hypersurface $D(h_1 + bh_2)$ does not
contain  any of the components of the intersection of the $D(f_i)$'s, $i \leq k$.
Let  $f_{k+1} = h_1 +bh_2$.
\end{proof}

\begin{cor}
Let $X$ be a projective variety over an infinite field $k$, $\phi:X \to X$
a morphism
 and ${{L}}$  an ample line bundle on $X$ such that 
$\phi^*({{{L}}}) \cong {{{L}}}^{\te d}$ for some $d \geq 1$. 
Then there exists an embedding $\iota$ of $X$ in $\bP^N_k$ and a morphism 
$\psi:\bP^N_k \to \bP^N_k$ such that  $ \psi \circ \iota  = \iota 
\circ \phi$.\label{cor:embed}
\end{cor} 

\begin{proof}
If $d = 1$ then $f$ must be an automorphism. Replacing ${{L}}$
by a very ample tensor power ${{L}}' = {{L}}^{\otimes n}$, one
sees that the natural embedding $\iota:X \to \bP(H^0(X, {{L}}))$
and the automorphism $\psi$ of $\bP(H^0(X, {{L}}'))$
induced by $f^*: H^0(X, {{L}}') \to H^0(X, {{L}}')$ have the required
properties.

If $d>1$, we use the  results of Mumford (\cite{mum-quad}, Theorems 1 and 3) to replace
${{L}}$ by a tensor power ${{L}}' = {{L}}^{\otimes n}$
so that the conditions  of Proposition \ref{prop:embed}
are satisfied by $X$, $\phi$ and ${{L}}'$.
\end{proof}

Our main application of Proposition \ref{prop:embed} is to show that the Boundedness Conjecture
of Morton and Silverman for $\bP^N$ implies the uniform boundedness of torsion for abelian 
varieties, using an idea of Bjorn Poonen.
Recall that for a map of sets $f:X \to X$, a point $x \in X$ is called \emph{preperiodic} for
$f$ if $f^n(x) = f^m(x)$ for some $m>n\geq 0$;
it is called \emph{periodic} if we can take $n=0$. 

\begin{conj}[Morton-Silverman  \cite{morton-silverman}] \label{conj:ms}
For all positive integers 
$D, N, d$ with $d>2$, there exists an integer $\kappa(D,N,d)$
such that for each number field
$k$ of degree $D$ over $\Q$, and each morphism $\psi:{\bP}^N_k \to
{\bP}^N_k$ defined by homogenous polynomials of degree $d $ over $k$, 
the number of preperiodic points of $\psi$
in ${\bP}^N(k)$ is less than or equal to $\kappa(D,N,d)$.
\end{conj}

\vspace{3 mm}
For an abelian variety $A$ and $n \in \Z$ we denote by $[n]:A \to A$ the multiplication 
by $n$ map. If $n \neq  0$, then any point of $A(k)_{tors}$ is preperiodic for
$[n]$.

\begin{cor} \label{cor:1}
The Boundedness Conjecture implies that there exists a constant $\eta(D,g)$
such that for any abelian variety $A$ of dimension $g$ defined over a number
field $k$ of degree $D$ over $\Q$, $|A(k)_{tors}| \leq \eta(D,g)$.
\end{cor}
\begin{proof}
 Let $A$ be an abelian variety of dimension $g$ defined over a number
field $k$. By Zarhin's trick (\cite[Lemma 2.5]{zarhin-trick} or \cite[p.~205]{moret-bailly})
 $X = (A \times A^t)^4$ has a principal polarisation $\lambda$. By \cite[p.~121]{GIT2},
$2\lambda$ is represented by a canonical symmetric\footnote{The symmetry
follows from the construction in loc.~cit.} ample line bundle
${{M}}$ on $X$. Let ${{{L}}} = {{M}}^{\otimes 3}$; 
this is a symmetric and very ample line bundle \cite[p.~163]{av}. 
By a theorem of Kempf \cite{kempf}, 
the imbedding of $X$ in ${\bP}(H^0(X,{{{L}}}))$ is projectively normal and the 
homogenous ideal of the image of $X$ is 
generated by its elements of degree two and three. Furthermore, since ${{{L}}}$ 
is symmetric $[n]^*{{{L}}} \cong {{{L}}}^{\te n^2}$ for all $n \in \Z$. 

Let $n =2$, so $d=2^2 =4$. We apply
Proposition \ref{prop:embed} to $X $,  ${{{L}}}$ 
and $\phi = [2]$.
By the Riemann-Roch theorem for abelian varieties (\cite[p.~150]{av})
$h^0(X,{{{L}}}) = (8g)!\  6^{8g}$. It follows from the
Boundedness Conjecture  that 
$|X| = |(A \times A^t)^4(k)_{tors}| \leq \kappa(D, (8g)! \ 6^{8g}-1,4)$. Thus we may 
take  $\eta(D,g)= \kappa(D, (8g)!\ 6^{8g}-1,4)^{1/4}$.
\end{proof}

\begin{rem}
If we assume Conjecture \ref{conj:ms} only for periodic points
 then Corollary \ref{cor:1} still holds, but we get a slightly weaker bound:
A  point in $A(k)_{tors}$ is periodic for $[n]$ iff its order 
is prime to $n$.
Applying the method
 of proof of Corollary \ref{cor:1}  to the maps $[2]$ and $[3]$, 
we get uniform bounds for the 
torsion of order prime to $2$ and the torsion of order prime to $3$.
Then $|A(k)_{tors}|$ is  bounded by the product of these bounds.
\end{rem}

\begin{rem}
Let  $\psi:{\bP}^N_k \to {\bP}^N_k$ be a morphism. By the functorial property
of restriction of scalars, there is an induced morphism
$\psi': \operatorname{Res}^{k}_{\mb{Q}}(\mb{P}^N_k) \to \operatorname{Res}^{k}_{\mb{Q}}(\mb{P}^N_k)$.
Preperiodic points for $\psi$ in ${\bP}^N(k)$ then correspond bijectively
to preperiodic points for $\psi'$ in $\operatorname{Res}^{k}_{\mb{Q}}(\bP^N_k)(\mb{Q})$.
One can apply
Proposition \ref{prop:embed} to $\psi'$ and the line bundle ${{L}}$
on $\operatorname{Res}^{k}_{\mb{Q}}(\mb{P}^N_k)$ induced by $\mc{O}(1)$,
 uniformly as $\psi$ varies over
all morphisms ${\bP}^N_k \to {\bP}^N_k$ with $N, D, d$ fixed and $[k:\mb{Q}] = D$,
to see that the special case of Conjecture \ref{conj:ms} for $D= 1$ and all $d,N$, 
in fact implies the general case.
\end{rem}

\section{} \label{sec:2}

We now consider some uniformity properties of 
the images of the set of  \emph{all} rational points
for varying morphisms of projective spaces. For the rest of this section
$k$ will be a number field and  $d>1, N$ will be positive integers. We shall also
abuse convention and say that a self map of $\bP^N_k$ is of degree
$d$ if it is defined by homogenous polynomials of degree $d$ over $k$.

Given a morphism
$f:{\mb{P}}^N_k \to {\mb{P}}^N_k$
of degree $d$ over $k$,
we may form the inverse system of sets
\[
 \cdots  \stackrel{f}{\to} \mb{P}^N(k) \stackrel{f}{\to} \cdots \stackrel{f}{\to} \mb{P}^N(k) \stackrel{f}{\to} \mb{P}^N(k) .
\]
Using the theory of heights one can show that the set 
$\varprojlim \  \mb{P}^N(k)$ is finite, and its 
points correspond in a natural way with the periodic points of $f$.
It is then natural to ask the following 
\begin{ques}
Let $\{f_i\}_{i=1}^{\infty}$ be a sequence of morphisms
$f_i:{\mb{P}}^N_k \to {\mb{P}}^N_k$
of degree $d$ over $k$. Consider the inverse system of sets
\[
 \cdots  \stackrel{f_n}{\to} \mb{P}^N(k) \stackrel{f_{n-1}}{\to} \cdots \stackrel{f_2}{\to} \mb{P}^N(k) \stackrel{f_1}{\to} \mb{P}^N(k) .
\]
Is the set $\varprojlim \  \mb{P}^N(k)$ finite? \label{ques:1}
\end{ques}

If true, the following  would have a positive answer.
\begin{ques}
Does there exist a finite subset $S(d,N,k)$ of ${\mb{P}}^N(k)$
such that for all morphisms $f:{\mb{P}}^N_k \to {\mb{P}}^N_k$ 
of degree $d$ over
$k$, $S(d,N,k)$ is not contained in $f({\mb{P}}^N(k))$? \label{ques:2}
\end{ques}

\begin{prop} \label{prop:lang}
Lang's conjecture on rational points on varieties of general type
implies that Question \ref{ques:2} has a positive answer.
\end{prop}

\begin{proof}
We shall consider the cases $N=1$ and $N>1$ separately.

\vspace{2mm}
\underline{$N=1$}. 
$PGL(2,k)$ acts $3$-transitively on $\bP^1(k)$, so for a given rational 
function $f$ we may find an automorphism $h$ of $\bP^1_k$ such that
$hf:\bP^1_k \to \bP^1_k$ is not ramified over $0$ or $\infty$. Since
rational functions $f$ of a fixed degree $d$ 
 depend on a finite number of parameters, we may
find a finite set of automorphisms $\{h_i\}$ such that for any rational
function $f$ of degree $d$, at least one of the rational functions
$h_if$ is not ramified over $0$ or $\infty$.

Consider the smooth projective model $X$ of the curve defined by the
equation $y^3 = f(x)$, where $f$ is a rational function of degree $d > 1$ 
not ramified over $0$ or $\infty$. $X$ is ramified over all the zeros and poles
of $f$ and the ramification degree at such points is $3$. Applying the
Riemann-Hurwitz formula, one sees that the genus of $X$ is at least $2$.
The genera of all curves obtained in this way (for a fixed $d$)
is also bounded above, since
there are only finitely many possibilities for the ramification type.
By the theorem of Caporaso-Harris-Mazur \cite{chm}, Lang's conjecture
implies that the number of $k$-rational points on all such curves
is bounded by a constant $C_k$. Let $r_d$ be the maximum number of points
of $\bP^1$ over which a rational function of degree $d$ can be ramified,
and let $T$ be a subset of cubes in $\Q^*$ with $|T| > C_k + r_d$.
Let $S = \cup_i h^{-1}_i(T)$. If $S \subset f(\bP^1(k))$ for some $f$
of degree $d$, the above considerations would imply that one of the
curves constructed above has more than $C_k$ $k$-rational points,
a contradiction.

\vspace{2mm}
\underline{$N > 1$}. For a fixed $f$,  Bertini's theorem implies that $f^{-1}(L)$
is smooth and irreducible for a general line $L$ in $\mb{P}^N$. It follows that
for a fixed degree $d$, there exists a finite set of lines $\{L_i\}_{i=1}^{r(d)}$
such that for any $f$ of degree $d$ at least one of the curves $\{f^{-1}(L_i)\}_{i=1}^{r(d)}$
is smooth and irreducible.

Suppose $f$ is general. By the adjunction formula the degree of the branch
locus of $f$ is $(N+1)(d -1)$, so the discriminant locus is 
an irreducible 
(since $f$ is general) hypersurface of degree $d^{N-1}(N+1)(d -1)$.
It follows that $f|_{f^{-1}(L_i)}:f^{-1}(L_i) \to L_i$ is of degree $d^N$ and branched
over $d^{N-1}(N+1)(d -1)$ points of $L_i$. From the Riemann-Hurwitz
formula we see that  the genus of $f^{-1}(L_i)$ is at least $2$ if
\[
d^{N-1}(N+1)(d -1) \geq 2d^N + 2 .
\]
This holds for $d \geq 4$ if $N=2$, $d \geq 3$ if $N=3$ and for 
 $d \geq 2$ if $N \geq 4$. Since the genus remains constant in families,
it follows that this will hold for arbitrary $f$ as long as $f^{-1}(L_i)$
is smooth.

For $d,N$ satisfying the above conditions, we can
apply the Caporaso-Harris-Mazur theorem to obtain subsets $S_i$ of 
$L_i(k)$, $i = 1, \ldots,r(d)$, such that $\bigcup S_i$  is not contained
in $f(\mb{P}^N(k))$ for any $f$ of degree $d$. For the exceptional pairs,
we repeat the above arguments using,instead of lines, smooth degree $2$ curves which are rational
over $k$.

\end{proof}

It would be very interesting to find a complete proof i.e.~not relying on any conjectures. 

Question \ref{ques:2} is evidently equivalent to the following which we state in order
to motivate the succeeding questions.

\begin{ques}
Let $h:\mb{P}^N(k) \to \R$ be a logarithmic Weil height.
Does there exist a constant $C(d,N,k) \in \mb{R}$ such
that $\{x \in \mb{P}^N(k) \ | \ h(x) < C(d,N,k)\}$ is not contained
in $f(\mb{P}^N(k))$ for all $f$ as above? \label{ques:3}
\end{ques}
One could also ask for the constant to depend only on the degree of $k$ over $\Q$.

In fact,
some computer calculations for 
the case $k=\Q$, $N=1$ and morphisms of the form $f(x) = x^2 + c, c \in \Q$
suggest that the following might be true.

\begin{ques}
For a morphism $f: \mb{P}^N_k \to \mb{P}^N_k$ 
of degree $d$ over $k$ and $c \in \mb{R}$,
let $N(f,c) = \#\{ x \in \mb{P}^N(k) \ | \ h(f(x)) \leq c \ \}$ and let
$N(c) = \underset{\deg(f)=d}{sup} \ N(f,c)$. 
Let $M(c) = \#\{ x \in \mb{P}^N(k) \ |\  h(x) \leq c\}$. Then is
\[
\lim_{c \to \infty} \frac{log(M(c))}{log(N(c))} = d \ ?
\]
\end{ques}

A generalisation of the Morton-Silverman Boundedness Conjecture (Conjecture \ref{conj:ms})
is the following:

\begin{ques}
For a morphism $f: \mb{P}^N_k \to \mb{P}^N_k$
of degree $d$ over $k$, 
let $h_f:\mb{P}^N(k) \to \R $ denote the
corresponding canonical height \cite[p.~30]{serre-mw}.
For $c \in \mb{R}$, let $R(c) = \underset{\deg(f)=d}{sup} \ \#\{ x \in \mb{P}^N(k) \ | \ h_f(x) \leq c\ \} $. Then is $R(c) < \infty$ for all $c$? If so, then with $M(c)$ as above, is
\[
\lim_{c \to \infty} \ \frac{log(R(c))}{log(M(c))} = 1 \ ?
\]
\end{ques}
It would be desirable to understand the precise relation between the previous two
questions. For example, does either of the questions imply the other?

If the previous question has a positive answer, then 
using the methods of Section \ref{sec:1}  one would get 
a positive answer to the following:
\begin{ques}
For each $c \in \mathbb{R}$,
does there exist  an integer $C(k,g,D, c)$
such that for any 
abelian variety $A$ of dimension $g$ over $k$ and an ample symmetric
line bundle ${{L}}$ of degree $D$ on $A$,
\[
\#\{x \in A(k) \ | \ h_{{{L}}}(x) \leq c \ \} < C(k,g,D, c) \ ?
\]
Here $h_{{{L}}}$ denotes the canonical height on $A(k)$ associated to the
line bundle ${{L}}$.
\end{ques}

\section{} \label{sec:3}

The existence of a self map of a variety satisfying the conditions of 
Corollary \ref{cor:embed} with $d>1$ imposes strong restrictions on 
the geometry of such varieties. The main result of this section,
Theorem \ref{thm:kodzero},
gives a complete classification in the case of smooth varieties of non-negative  Kodaira dimension.
We then state some questions concerning the general case.


\begin{lem}
Let $X$, ${{{L}}}$ and $\phi$ be as in Corollary \ref{cor:embed} with $d>1$. 
Assume further that
$\phi$ is separable and the Kodaira dimension of $X$ is $\geq 0$ i.e. 
$H^0(X,\omx^{\te n}) \neq 0$ for some $n>0$,
where $\omx$ is the canonical line bundle of $X$. Then $\phi$ is \'etale\footnote{This is probably
well-known.}, 
the class of $\omx$ is
torsion in $Pic(X)$ and $X$ does not contain any positive dimensional subvarieties with 
finite algebraic fundamental group.
\end{lem}
\begin{proof}
Let $n >0$ be such that $H^0(X,\omx^{\te n}) \neq 0$. Since  $\phi$ is
separable, $\omx = \phi^*{\omx} \te \mc{O}_X(R)$ where $R$ is the ramification
divisor.
Since $R$ is an effective divisor, it follows that any section of  
$H^0(X, \omx^{\te n})$ is in the image of $H^0(X, \omx^{\te n} \te \mc{O}_X(-nR))$. 
By iterating $\phi$ it follows that
\[
H^0(X, \omx^{\te n}) \subset \bigcap_{i\geq 0} H^0(X, \omx^{\te n} \te \mc{O}_X(- n(R + \phi^*(R)
 + \cdots + (\phi^{i-1})^* (R)))).
\]
The right hand side is zero unless $R=0$, so
$\phi$ must be \'etale.

Let $P({{{L}}},x)$ denote the Hilbert polynomial of ${{{L}}}$. We apply the
Grothendieck-Riemann-Roch 
theorem (see for example \cite[Theorem 15.2] {fulton-it})
to the morphism $\phi$. Since  $\phi$ is \'etale, $\phi^*(T_X) = T_X$,
so the relative tangent bundle is trivial.
Thus, for any integer $x$ we get
\begin{equation}
ch(\phi_*({{L}}^{\te dx})) = \phi_*(ch({{L}}^{\te dx})) . \label{eq-grr}
\end{equation}
Since 
$\phi^*({{L}}) \cong {{L}}^{\te d}$, the projection formula implies that
$\phi_*({{L}}^{\te dx}) \cong {{L}}^{\te x} \te \phi_*(\mc{O}_X)$.
The Hirzebruch-Riemann-Roch theorem (\cite[Corollary 15.2.1]{fulton-it})
then gives:
\begin{multline*}
  P({{{L}}}, dx) = \int_X ch({{L}}^{\te dx})\cdot td(T_X)  =
 \int_X \phi_*(ch({{L}}^{\te dx})\cdot td(T_X)) = \\
 \int_X \phi_*(ch({{L}}^{\te dx}))\cdot td(T_X) = 
 \int_X ch(\phi_*({{L}}^{\te dx}))\cdot td(T_X) = \\
 \int_X ch({{L}}^{\te x}) \cdot ch(\phi_*(\mc{O}_X)) \cdot td(T_X) =
 d^g \int_X ch({{L}}^{\te x})\cdot td(T_X) = 
 d^g P({{{L}}}, x)
\end{multline*}
with $g = \dim(X)$.
Here the third equality follows from the projection formula since $\phi^*(T_X) = T_X$,
the fourth from \eqref{eq-grr}, the fifth from multiplicativity of $ch$,
the sixth from $ch(\phi_*(\mc{O}_X)) = deg(\phi) = d^{\dim(X)}$ in $A^*(X)_{\Q}$,
which holds since  $\phi$ is \'etale (and the first and last follow from 
HRR).

Since $P({{{L}}},x)$ is a polynomial of degree $g$, it follows that
 $P({{{L}}},x) = P({{{L}}},1)\, x^g$. 
Examining the terms in the Hirzebruch-Riemann-Roch formula
we see that $\int_X (c_1({{{L}}})^{g-1} \cdot c_1(\omx)$, which is the
coefficient of $x^{g-1}$ in $P({{{L}}},x)$, must be $0$. Since 
${{{L}}}$ is ample and some multiple of $\omx$ has a section, this can
only happen if the zero set of the section is empty. Thus  $\omx$ is
torsion.

Let $Y \subset X$ be a subvariety of dimension $e$, $0<e \leq g$, with finite algebraic
fundamental group, say of order $f$.
Let $c = \int_X(c_1({{{L}}})^e \cdot Y) > 0$ and choose $m>0$
so that $cf/d^{e m} < 1$. Let $Y'$ be an irreducible component of
$(\phi^m)^{-1}(Y)$. It follows that $deg \ \phi^m |_{Y'} \leq f$.
 By the projection formula 
\[ 
d^{e \cdot m} \int_X  (c_1({{{L}}})^e \cdot Y') = \int_X 
(c_1({{{L}}}^{\te d^m})^e \cdot Y') \leq  \int_X (c_1({{{L}}})^e \cdot Y) 
\]
Thus  $\int_X (c_1({{{L}}})^{e} \cdot Y') \leq c f/d^{e m} < 1$. This
is a contradiction because ${{{L}}}$ is ample so the left hand side must be 
a positive integer.
\end{proof}

The classification result that we obtain is the following

\begin{thm} \label{thm:kodzero}
Let  $X$ be a smooth projective variety of non-negative Kodaira dimension, 
${{{L}}}$ an ample
line bundle on $X$  and $\phi$ a self map of $X$ such that
$\phi^*({{L}}) \cong {{L}}^{\te d}$ for some $d>1$, all over an
algebraically closed field $k$ of characteristic zero\footnote{It seems likely
that this also holds in positive characteristic if we assume that the map $\phi$ is
separable.}.
Then $X$ is isomorphic to a quotient $A/G$
of an abelian variety $A$ by a finite group $G$ acting fixed point freely. Furthermore,
all such quotients, 
over an algebraically closed field of arbitrary characteristic,
have line bundles ${{L}}$ and self maps $\phi$ as above.
\end{thm}
\begin{proof}
By the previous lemma we have that $\omx$ is torsion. By a theorem of
Beauville \cite[Th\'eor\`me 2]{beauville-chern} 
it follows that there exists an abelian variety $A$, a simply 
connected smooth projective variety $Z$ and a finite group $G$ acting
on $A \times Z$ without fixed points such that $X \cong (A \times Z)/G$.
We need to show that $Z$ is a point. Since $Z$ is simply connected, there
are no nonconstant morphisms from $Z$ to $A$. It follows that if $p_1$
denotes the projection from $A \times Z$ to the first factor, then
$p_1(g(\{a\} \times Z))$ is a point for all $g \in G$. Thus $G$ acts on
$A$ (not necessarily faithfully) and we get a surjective morphism
$\bar{p_1}:(A \times Z)/G \to A/G$ such that the diagram
$$
\diagram
A \times Z \dto_{p_1} \rto & (A \times Z)/G \dto_{\bar{p_1}} \\
A \rto & A/G \\
\enddiagram
$$
commutes. For all $a \in A$, $\bar{p_1}^{-1}(a)$ is isomorphic to a quotient
of $Z$ by a  finite group and hence has a finite fundamental group. Lemma 1
then implies that $Z$ must be a point.

Now suppose that $X \cong A/G$, with $A$, $G$ as above. Let $V = A(k)/A(k)_{tors}$; this is
a vector space over $\Q$ and the action of $G$ on $A$ induces an action of $G$ on
$V$ which factors through the semi-direct product of $GL(V)$ and $V$. 
Since $G$ is finite, this action must have
a fixed point: if $\bar{x} \in V$ is any point, the point
 $|G|^{-1}\sum_{g \in G} g \bar{x}$ is fixed by $G$.
By translation,  we may assume that $\bar{0} \in V$
is a fixed point
of the action of $G$ on $V$ i.e.~for all $g \in G$, $g(0) \in A(k)_{tors}$.
Let $n$ be an integer such that $n \equiv 1 \mod ord(g(0))$ for all $g \in G$.
Then the multiplication by $n$ map on $A$, $[n]$, commutes with the action of $G$, hence
descends to a morphism $\phi_n:X \to X$. Let ${{M}}_1$ be an ample, symmetric line bundle
on $A$. Let $m$ be the $g.c.d$.~of all $ord(g(0))$, $g \in G$, and let 
${{M}}  = \bigotimes_{g \in G} g^*({{M}}_1)^{\te m}$. Then $[n]^*({{M}}) \cong {{M}}^{\te n^2}$
and moreover $g^*({{M}}) = {{M}}$ for all $g \in G$. Since the
kernel of the pullback map $Pic(X) \to Pic(A)$ is finite, it follows
that a suitable positive 
tensor power of ${{M}}$  descends
to an ample line bundle ${{L}}$ on $X$ such that $\phi_n^*({{L}}) \cong {{L}}^{\te n^2}$.
\end{proof}

We thus have a reasonably complete description of varieties of non-negative Kodaira
dimension which have self maps and line bundles as above. However, the general case
seems much more difficult.
As a beginning,
one may ask the following:
\begin{ques}
Let $X$ be an $N$-dimensional smooth projective variety over an algebraically closed
field of characteristic zero. Suppose $Pic(X) \cong \Z$, $X$ has a 
a self map of degree $> 1$ and the anti-canonical bundle of $X$ is ample. 
Then is $X$ isomorphic to  $\bP^N$?
\end{ques} 

The question has a positive answer for
$N=1, 2$ as can be easily checked from the classification of smooth
projective curves and surfaces. It follows from results of Paranjape-Srinivas
\cite{kap-sri}
that a   variety of the form $G/P$, where $G$ is a simple linear algebraic group and
$P$ is a maximal parabolic subgroup (such varieties have $Pic(G/P) \cong \Z$)
has a self map of degree $> 1$ iff $G/P \cong \bP^N$ for some $N$. The same
holds for smooth projective hypersurfaces $H$ with $\dim(H)>2$ (these also
have $Pic(H) \cong \Z$) by
results of Beauville \cite{beau-hyp}.

One may ask the same question for separable morphisms of
varieties over algebraically closed fields of arbitrary characteristic. 
If one also allows singular varieties, there are other examples. But perhaps
all such, at least those which are normal, are toric varieties. (It is
well known that
projective toric varieties do have self maps and line bundles verifying the
conditions of Corollary \ref{cor:embed}).

The following is a possible generalisation. 
\begin{ques}
Let $X$ be a smooth projective rationally connected variety over an
algebraically closed field of characteristic zero.
Suppose $X$ has a self map $\phi$ and an ample line bundle ${L}$
such that
$\phi^*L \cong L^{\te d}$, for some
$d >1$. Then is $X$ a toric variety?
\end{ques}
It is not difficult to check this (using classification) if  $\dim(X) \leq 2$.
A more general result has been proved by N.~Nakayama \cite{nakayama-ruled}.

\vspace{2mm}
One may also try to classify \'etale maps of varieties of arbitrary Kodaira dimension.
\begin{ques}
Let $X$ be a smooth projective variety over an algebraically closed field and
$f:X \to X$ an \'etale map. Does there exist a finite \'etale cover $Y$
of $X$, a self map $g:Y \to Y$, a smooth projective morphism $\pi:Y \to S$ 
whose fibres are abelian varieties, and
an \emph{automorphism} $h:S \to S$ such that the following diagrams commute?
\[
\xymatrix{
Y \ar[r]^g \ar[d] & Y \ar[d] \\
X \ar[r]^f & X \\
}
\hspace{1in}
\xymatrix{
Y \ar[d]_{\pi} \ar[r]^g & Y \ar[d]^{\pi} \\
S \ar[r]^{h} & S\\
}
\]
\end{ques}
Again, this can be easily checked for curves and surfaces. 
Note that this is also true for $X$ of the from $A/G$ considered above:
for such $X$, among all \'etale maps $B \to X$ with
$B$ an abelian variety, there is a unique (upto isomorphism) map of minimal degree.
Some general results for $X$ of dimension $3$ have been obtained by 
Y.~Fujimoto \cite{fujimoto}.

\section{} \label{sec:4}

Let $X$ be a projective variety over an algebraically closed field $k$, $\phi:X \to X$
a morphism and ${{L}}$  an ample line bundle on $X$ such that 
$\phi^*({{{L}}}) \cong {{{L}}}^{\te d}$ for some $d > 1$. Using the Lefschetz trace
formula and a result of Serre \cite[Th\'eor\`eme 1]{serre-kahler}, one can  show that
if $X$ is smooth and $k$ is of characteristic zero then $\phi$ has infinitely many periodic
points in $X(k)$. However, we shall show below that under  more general conditions
the set of periodic points is always Zariski dense in $X$. The main ingredient
is a generalisation of the Lang-Weil estimates due to E.~Hrushovski \cite{hrush}.

\begin{thm} \label{thm:density}
Let $X$ be a projective variety over an algebraically closed field $k$, $\phi:X \to X$
a dominant morphism and  ${{L}}$ a line bundle on $X$ such that
$\phi^*{L} \otimes {L}^{- 1}$ is ample. Then the subset of $X(k)$
consisting of periodic
points of $\phi$ is Zariski dense in $X$. \label{thm:ZD}
\end{thm}

\begin{rem}
From the arguments below one may extract an elementary  proof of the
fact that under the conditions of the theorem the set of all preperiodic
points is Zariski dense---this does not need Hrushovski's theorem.
\end{rem}

\begin{cor}
Let $\phi:\bP^N \to \bP^N$ be a morphism of degree $d > 1$ over an algebraically closed
field $k$. Then the set of periodic points of $f$ is Zariski dense in $\bP^N$.
\end{cor}

\begin{rem}
Simple examples (e.g. $N=1$, $\phi(x) = x^2$) show that if $k= \mathbb{C}$, 
the set of periodic points need not
be dense in the analytic topology.
\end{rem}

Theorem \ref{thm:ZD} will be deduced from the following result, whose proof is 
due to Bjorn Poonen.

\begin{prop} \label{prop:density-finite}
Let $X$ be an 
algebraic  variety over $\mb{F}$, 
the algebraic closure of a finite field and 
$\phi:X \to X$ a finite surjective morphism. 
Then the subset of $X(\mb{F})$ consisting of
periodic points of $\phi$ is Zariski dense in $X$.
\end{prop}
\begin{proof}
Let $Y$
be the Zariski closure of the set of periodic points
of $\phi$ in $X(\mb{F})$ and suppose $Y \neq X$.
Let $q = p^n$, $p = char(\mathbb{F})$, be  such that 
 $X$ as well as $\phi$ are defined over the subfield
$\mathbb{F}_q$ of $\mathbb{F}$ consisting of $q$ elements. 
Let $\sigma$ denote the Frobenius morphism
of $X$ (which raises the coordinates of a point $x \in X$ to their $q$'th powers)
and let $\Gamma_{\phi}$ (resp. $\Gamma_m$) denote the
graph of $\phi$ (resp.~$\sigma^m$) in $X \times X$. 
Let $U$ be an irreducible affine open subset of $X -Y$ also defined over $\mathbb{F}_q$
and let $V = \Gamma_{\phi} \cap  (U \times U)$.
By a theorem of Hrushovski \cite{hrush}
(see \cite[Conjecture 2]{chat-hrush}\footnote{Z.~Chatzidakis has informed the author that
in the statement of \cite[Conjecture 2]{chat-hrush}, hence in Hrushovski's theorem,
 the projections only need to be
generically onto, not necessarily onto.} 
for the statement)
there exists $m >0$ such that $(V \cap \Gamma_m)(\mb{F}) \neq \emptyset$
i.e.~there exists $u \in U(\mb{F})$ such that $\phi(u) = \sigma^m(u)$.
Since $\phi$ is defined over $\mb{F}_q$, it follows that $u$ is a periodic
point of $\phi$. This contradicts the definition of $Y$ and $U$, so the
proof is complete.

\end{proof}

\begin{lem}
Let $X$ be a variety over a field $k$ and $\phi:X \to X$ a dominant projective
morphism. Then $\phi$ is a finite morphism.\label{lem:dom}
\end{lem}

\begin{proof}
We may assume that $k$ is algebraically closed.

Let $H_*(X)$ denote the \'etale homology of $X$
 with cofficients in $\Q_l$,
$l$ a prime with $(char\ k, l) =1$  (i.e.~the linear dual
of $H^*_{et}(X, \Q_l)$)
and let $P_*(X)$ denote the subspace
of $H_*(X)$ generated by classes of complete subvarieties of $X$. $\phi$ induces
a map $\phi_*:H_*(X) \to H_*(X)$ 
which preserves $P_*(X)$ i.e.~$\phi_*(P_*(X)) \subset P_*(X)$. Suppose $\phi$
is not finite, so there is an irreducible  positive dimensional subvariety $Y$
contained in
a fibre of $\phi$. The class of Y in $P_*(X)$, $[Y]$, is non-zero because $\phi$
is projective: $c_1({{L}})^{\dim(Y)} \cap [Y] \in H_0(X)$ is non-zero for
any relatively ample line bundle ${{L}}$ on $X$. However, $\phi_*([Y]) = 0$
because  $\dim(Y) > 0$ and $\dim(\phi(Y)) = 0$. Moreover,  $\phi_*:P_*(X) \to
P_*(X)$ is surjective since 
$\phi$ is dominant and proper, so any complete irreducible subvariety $Z$ of $X$ is
of the form $\phi(Z')$ for a complete subvariety of $X$ with $\dim(Z') = \dim(Z)$.
This is a contradiction because $H_*(X)$, hence $P_*(X)$,
 is a finite dimensional $\Q_l$ vector space.
\end{proof}

\begin{proof}[Proof of Theorem \ref{thm:ZD}.]
Since $X$ is projective and $\phi$ is dominant, it follows from Lemma
\ref{lem:dom}
that $\phi$ is a finite
morphism.
By Proposition \ref{prop:density-finite}, 
the theorem holds if $k$ is the algebraic
closure of a finite field. For a  general $k$, $X$, $\phi$ and ${{L}}$
will be  defined over a finitely
generated subfield $k'$ of $k$, so we may assume that the transcendence degree of $k$
over its prime field is finite.  It is easy to see that we may also
assume that $X$ is irreducible.


Now consider the following situation: $\mathbf{X}$ is a projective
scheme, flat  over a d.v.r.~$R$ with algebraically closed residue field $K$,  
${\Phi}$ is a dominant selfmap of $\mf{X}$ over $R$ and $\mf{L}$
a  line bundle  on $\mf{X}$ such that
$\Phi^*{\mf{L}} \te \mf{L}^{-1}$ is relatively ample. 
Let $X$ be the special fibre of $\mf{X}$ and $\phi$ the
restriction of $\Phi$ to $X$.
The set of periodic points of $\phi$
of period dividing a positive integer $n$  can be viewed as
the set of points in $(\Delta_X \cap \Gamma_{\phi^n})(K)$, where $\Delta_X$
is the diagonal and $\Gamma_{\phi^n}$ is the graph of $\phi^n$ in $X \times X$.
The hypothesis on the  line bundle $\mf{L}$
implies that this set is always finite: The line bundle $p_1^*{{L}}^{-1} \te p_2^*{{L}}$
on $X \times X$,
where ${{L}} = \mf{L}|_X$,
restricts to a line bundle on $\Delta_X \cap \Gamma_{\phi^n}$ which is
both ample and trivial.

Suppose that $X$ is reduced and the set of periodic points of $\phi$ is 
Zariski dense in $X$, so the set of periodic points of
$\phi$ which are smooth points on $X$  is also Zariski dense
in $X$. For a positive integer $n$, consider the subscheme
$\Delta_{\mf{X}} \cap \Gamma_{\Phi^n}$ of $X \times_R X$, where
$\Delta_{\mf{X}}$ is the diagonal and $\Gamma_{\Phi^n}$ is the 
graph of $\Phi$ in $X \times_R X$. If $x \in X$ is a periodic
point of $\phi$ which is also a smooth point on $X$, it follows
from \cite[Theorem 3, p.110]{serre-loc} that $(x,x)$ is contained in
a closed subscheme of $\Delta_{\mf{X}} \cap \Gamma_{\Phi^n}$ of dimension
at least 1. Since $(\Delta_X \cap \Gamma_{\phi^n})(K)$ is finite,
it follows that this subscheme must intersect the generic fibre of $X \times_R X$. So
$x$ can be lifted to a periodic point of $\Phi$ on any geometric generic fibre
of $X_{/R}$. It follows that the set of periodic points of $\Phi$ on any
geometric generic fibre are also Zariski dense.

The proof is completed by induction on the transcendence degree of
$k$, using models over discrete valuation rings as above to increase the transcendence degree by
one at each step. For $k$ of characteristic zero, we also
need to use one such  model to go from the algebraic closure of a finite
field to $\overline{\mb{Q}}$.


\end{proof}

\vspace{2mm}
\emph{Acknowledgements.}
The results of Section \ref{sec:1}
were motivated by a conversation with Bjorn Poonen who  observed that the
Morton-Silverman conjecture for degree $4$ maps of $\bP^1$ implies the
uniform boundedness of torsion for elliptic curves, and asked if that could
be generalised. I also thank him for informing me of Hrushovski's theorem
\cite{hrush}, its application to the proof of Proposition \ref{prop:density-finite},
and for his permission to include it here. This note has  benefited from
his numerous comments on a preliminary version, in particular the first part of
the proof of Proposition \ref{prop:lang} is a simpler variant, due to him,
of my original proof. Finally, I  thank Zo\'e Chatzidakis for some clarifications concerning
\cite[Conjecture 2]{chat-hrush}.


\end{document}